\documentclass[a4paper,11pt]{article}
\usepackage[latin1]{inputenc}
\usepackage{amsmath}
\usepackage{amsfonts}
\usepackage{amssymb}
\usepackage{amscd}
\usepackage{dsfont}
\usepackage{mathrsfs}
\usepackage{graphicx}
\usepackage{setspace}
\usepackage{amsthm}
\usepackage[all]{xy}

\newtheorem*{theorem*}{Theorem}
\newtheorem*{proposition*}{Proposition}
\newtheorem{theorem}{Theorem}[section]

\newtheorem{lemma}[theorem]{Lemma}
\newtheorem{proposition}[theorem]{Proposition}

\newtheorem{remark}[theorem]{Remark}

\newcommand{\disp}{\displaystyle}
\newcommand{\erre}{\mathds{R}}

\newcommand{\esse}{\mathds{S}}

\newcommand{\di}{\mathrm{d}} 

\newcommand{\ra}{\rightarrow}

\newcommand{\tc}{\,:\,}                                  

\begin{document}

\author{Debora Impera \and Luciano Mari \and Marco Rigoli}
\title{\textbf{Some geometric properties of hypersurfaces with constant $r$-mean curvature in Euclidean space}}
\date{}
\maketitle
\scriptsize \begin{center} Dipartimento di Matematica,
Universit\`a
degli studi di Milano,\\
Via Saldini 50, I-20133 Milano (Italy)\\
E-mail addresses: \\
debora.impera@unimi.it, \ luciano.mari@unimi.it, \
marco.rigoli@unimi.it
\end{center}

\begin{abstract}
Let $f:M\ra \erre^{m+1}$ be an isometrically immersed
hypersurface. In this paper, we exploit recent results due to the
authors in \cite{bimari} to analyze the stability of the
differential operator $L_r$ associated with the $r$-th Newton
tensor of $f$. This appears in the Jacobi operator for the
variational problem of minimizing the $r$-mean curvature $H_r$.
Two natural applications are found. The first one ensures that,
under a mild condition on the integral of $H_r$ over geodesic
spheres, the Gauss map meets each equator of $\esse^m$ infinitely
many times. The second one deals with hypersurfaces with zero
$(r+1)$-mean curvature. Under similar growth assumptions, we prove
that the affine tangent spaces $f_*T_pM$, $p\in M$, fill the whole
$\erre^{m+1}$.
\end{abstract}

 \maketitle

\normalsize

\section{Introduction}
In what follows $f \tc M^m \ra \erre^{m+1}$ will always denote a
connected, orientable, complete, non compact hypersurface of
Euclidean space. We fix an origin $o\in M$ and let
$r(x)=\mathrm{dist}(x,o)$, $x\in M$. We set $B_r$ and $\partial
B_r$ for, respectively, the geodesic ball and the geodesic sphere
centered at $o$ with radius $r$. Moreover, let $\nu$ be the
spherical Gauss map and denote with $A$ both the second
fundamental form and the shape operator in the orientation of
$\nu$. Associated with $A$ we have the principal curvatures
$k_1,\ldots, k_m$ and the set of symmetric functions $S_j$:
$$
S_j = \sum_{i_1<i_2<\ldots <i_j} k_{i_1}\cdot
k_{i_2}\cdot\ldots\cdot k_{i_j}, \qquad j\in \{1,\ldots, m\},
\quad S_0=1.
$$
The $j$-mean curvature of $f$ is defined
$$
H_0=1,\qquad {n \choose j} H_j = S_j,
$$
so that, for instance, $H_1$ is the mean curvature and $H_m$ is
the Gauss-Kronecker curvature of the hypersurface. Note that, when
changing the orientation $\nu$, the odd curvatures change sign,
while the sign of the even curvatures is an invariant of the
immersion. By Gauss equations and flatness of $\erre^{m+1}$ it is
easy to see that
$$
H_2= {m\choose 2}^{-1}S_2 = \frac 12 {m\choose 2}^{-1}
\mathrm{scal},
$$
where $\mathrm{scal}$ is the scalar curvature of $M$. The $j$-mean
curvatures satisfy the so-called Newton inequalities
$$
H_j^2\ge H_{j-1}H_{j+1},
$$
equality holding if and only if $p$ is an umbilical point (see
\cite{hardylittlpolya}). We stress that no restriction is made on
the sign of the $H_i$'s.
\begin{theorem}\label{main}
Let $f:M\ra \erre^{m+1}$ be a hypersurface such that, for some
$j\in\{0,m-2\}$, $H_{j+1}$ is a non-zero constant. If $j\ge 1$,
assume that there exists a point $p\in M$ at which the second
fundamental form is definite. Set
\begin{equation}\label{definV}
v_j(r)=\int_{\partial B_r} H_j, \qquad v_1(r)= \int_{\partial B_r}
H_1.
\end{equation}
where integration is with respect to the $(m-1)$-dimensional
Hausdorff measure of $\partial B_r$. Fix an equator $E\subset
\esse^m$ and suppose that either
\begin{equation}\label{assunz}
\begin{array}{ll}
(i) & \disp \int^{+\infty}\dfrac{\di r}{v_j(r)} = +\infty \qquad
\text{and}
\qquad H_1 \not\in L^1(M) \quad \text{or} \\[0.5cm]
(ii) & \disp\int^{+\infty}\dfrac{\di r}{v_j(r)} <+\infty \qquad
\text{and}
\\[0,4cm]
& \displaystyle \limsup_{r\ra
+\infty}\sqrt{v_1(r)v_j(r)}\int_r^{+\infty}\frac{\di s}{v_j(s)}
>\frac 12 \left[(j+1){m+1\choose j+2}H_{j+1}\right]^{-1/2}.
\end{array}
\end{equation}
Then, there exists a divergent sequence $\{x_k\} \subset M$ such
that $\nu(x_k)\in E$, where $\nu$ is the spherical Gauss map.
\end{theorem}
\begin{remark}
\emph{Up to changing the orientation of $M$, we can suppose that
the second fundamental form at $p$ is positive definite. As we
will see later in more detail, this has the remarkable consequence
that each $H_i$, $1\le i\le n$, is strictly positive at every
point of $M$. In particular, $v_1$ and $v_j$ are both strictly
positive and the requirements in \eqref{assunz} are meaningful.}
\end{remark}
\begin{remark}\label{remmain}
\emph{When $j=1$, the existence of an elliptic point $p\in M$ can
be replaced by requiring $H_2$ to be a positive constant, see
\cite{elbert} for details. The case $j=0$ has been considered in
\cite{bimari}.}
\end{remark}
We clarify the role of $(i)$ and $(ii)$ with some examples. First,
we deal with the case $j\neq 1$, and we assume that $v_j$ is of
order $r^k$ (resp $e^{kr}$), for some $k>0$. Then assumption
$(ii)$ requires that $v_1(r)$ is of order at least $r^{k-2}$ (resp
$e^{kr}$). Roughly speaking, $v_1$ has to be big enough with
respect to the other integral curvature $v_j$. Under additional
requirements on the intrinsic curvatures of $M$, standard volume
comparisons allow to control the volume of $\partial B_r$ and
$(ii)$ can be read as $H_1$ not decaying too fast at infinity.
When $j=1$, things are somewhat different. Indeed, $(ii)$ implies
that $v_1(r)$ does not \emph{grow} too fast, that is, loosely
speaking, it has at most exponential growth. This shows that two
opposite effects balances in condition $(ii)$. The same happens
for $(i)$ with $j=1$, as a consequence of Cauchy-Schwartz
inequality and coarea formula
$$
\left(\int^r_R \frac{\di
s}{v_1(s)}\right)\left(\int_{B_r\backslash B_R} H_1\right) \ge
(r-R)^2.
$$
Finally, we stress that $(i)$ and $(ii)$ are mild hypotheses as
they only involve the integral of extrinsic curvatures. In other
words, no pointwise control is required.\\
\par
Up to identifying the image of the tangent space at $p\in M$ with
an affine hyperplane of $\erre^{m+1}$ in the standard way, we can
also prove the following result:
\begin{theorem}\label{main2}
Let $f:M\ra \erre^{m+1}$ be a hypersurface with $H_{j+1}\equiv 0$.
If $j\ge 1$, assume $\mathrm{rank}(A)>j$ at every point. Define
$v_1,v_j$ as in \eqref{definV}. Then, under assumptions
\eqref{assunz} $(i)$ or $(ii)$, for every compact set $K\subset M$
we have
$$
\bigcup_{p\in M\backslash K}T_pM\equiv \erre^{m+1},
$$
that is, the tangent envelope of $M\backslash K$ coincides with
$\erre^{m+1}$.
\end{theorem}
\begin{remark} \emph{As we will see later, condition
$\mathrm{rank}(A)>j$ implies that $H_i>0$ for every $1\le i\le
j$.}
\end{remark}

\section{Preliminaries}
We start recalling the definition and some properties of the
Newton tensors $P_j$, $j\in\{0,\ldots, m\}$. They are inductively
defined by
$$
P_0= I, \qquad P_{j} = S_jI-AP_{j-1}.
$$
For future use, we state the following algebraic lemma. For a
proof, see \cite{barbosacolares}.
\begin{lemma}\label{algebrico}
Let $\{e_i\}$ be the principal directions associated with $A$,
$Ae_i=k_ie_i$, and let $S_j(A_i)$ be the $j$-th symmetric function
of $A$ restricted to the $(m-1)$-dimensional space $e_i^\perp$.
Then, for each $1\le j\le m-1$,
$$
\begin{array}{ll}
(1) & AP_j=P_jA; \\[0.2cm]
(2) & P_j e_i= S_j(A_i) e_i; \\[0.2cm]
(3) & \mathrm{Tr}(P_j) = \sum_i S_j(A_i) = (m-j)S_j; \\[0.2cm]
(4) & \mathrm{Tr}(AP_j) = \sum_i k_i S_j(A_i)= (j+1)S_{j+1}; \\[0.2cm]
(5) & \mathrm{Tr}(A^2P_j) = \sum_i k_i^2S_j(A_i) =
S_1S_{j+1}-(j+2)S_{j+2}.
\end{array}
$$
\end{lemma}
It follows from $(2)$ in the above lemma, and from the definition
of $P_m$ that $P_m=0$. Related to the $j$-th Newton tensor there
is a well defined, symmetric differential operator acting on
$C^\infty_c(M)$:
\begin{equation}\label{diver}
 L_ju =
\mathrm{Tr}(P_j\mathrm{Hess}(u))= \mathrm{div}(P_j\nabla u) \qquad
\forall \ u\in C^\infty_c(M),
\end{equation}
where the last equality is due to the fact that $A$ is a Codazzi
tensor in $\erre^{m+1}$, see \cite{chengyau}, \cite{ros}. $L_j$
naturally appears when looking for stationary points of the
curvature integral
$$
\mathcal{A}_j(M) = \int_M S_j \di V_M,
$$
for compactly supported volume preserving variations. These
functionals can be viewed as a generalization of the volume
functional. In fact, in \cite{barbosacolares} and \cite{elbert}
the stationary points of $\mathcal{A}_j$ are characterized as
those immersions having constant $S_{j+1}$. In the above mentioned
paper \cite{elbert}, M.F. Elbert computes the second variation of
$\mathcal{A}_j$ in more general ambient spaces and obtains in the
Euclidean setting the expression
$$
T_j = L_j + \big(S_1S_{j+1}-(j+2)S_{j+2}\big)
$$
for the Jacobi operator. In what follows we are interested in the
case of $L_j$ elliptic. There are a number of different results
giving sufficient conditions to guarantee this fact, and the next
two fit the situation of our main theorems.
\begin{proposition}
Let $M$ be an $m$-dimensional connected, orientable hypersurface
of some space form $N$. Then, $L_i$ is elliptic for every $1\le
i\le j$ in each of the following cases:
\begin{itemize}
\item[(i)] $M$ contains an elliptic point, that is, a point
$p\in M$ at which $A$ is definite (positive or negative), and
$S_{j+1}\neq 0$ at every point of $M$. Note that, up to changing
the orientation of $M$, we can assume $A_p$ to be positive
definite, and by continuity $S_{j+1}>0$ on $M$.
\item[(ii)] $S_{j+1}\equiv 0$ and $\mathrm{rank}(A)>j$ at every point of $M$.
\end{itemize}
Moreover, in both cases, every $i$-mean curvature $H_i$ is
strictly positive on $M$, for $1\le i\le j$.
\end{proposition}
For a proof of $(i)$ see \cite{barbosacolares}, while for $(ii)$ see \cite{houleite}.

From the above proposition, the requirements on $p$ and
$\mathrm{rank}(A)$ in the main theorems ensure ellipticity. As
stressed in Remark \ref{remmain}, when $j=2$ in \cite{elbert} it
is shown that the sole requirement $H_2>0$ implies the ellipticity
of $L_1$. In the assumptions of the above proposition, we can
define the $j$-volume of some measurable subset $K\subset M$ as
the integral
$$
\mathcal{A}_j(K) = \int_KS_j\di V_M.
$$
Hereafter, we restrict to the case $L_j$ elliptic. Given the
relatively compact domain $\Omega\subset M$, $L_j$ is bounded from
below on $C^\infty_c(\Omega)$ and, by Rellich theorem, for a
sufficiently large $\lambda$, $(L_j-\lambda)$ is invertible with
compact resolvent. By standard spectral theory, $L_j$ is therefore
essentially self-adjoint on $C^\infty_c(\Omega)$ (Theorem 3.3.2 in
\cite{vanishing}). Essential self-adjointness implies that
$C^\infty_c(\Omega)$ and $\mathrm{Lip}_0(\Omega)$ are cores for
the quadratic form associated to $L_j$. The first eigenvalue
$\lambda_1^{T_j}(\Omega)$, with Dirichlet boundary condition, is
therefore defined by the Rayleigh characterization
$$
\lambda_1^{T_j}(\Omega) = \inf_{ \scriptsize{\begin{array}{c} \phi\in \mathrm{Lip}_0(\Omega) \\[0.1cm] \phi \neq 0 \end{array} } }
 \frac{\int_\Omega \langle P_j(\nabla \phi),\nabla \phi\rangle - \int_\Omega (S_1S_{j+1}-(j+2)S_{j+2})\phi^2}{\int_\Omega
 \phi^2},
$$
where $\mathrm{Lip}_0(\Omega)$ can be replaced with
$C_c^\infty(\Omega)$. By the monotonicity property of eigenvalues
(or, in other words, since $L_j$ satisfies the unique continuation
property, \cite{aronsz}), if $\Omega_1$ is a domain with compact
closure in $\Omega_2$, and $\Omega_2\backslash\Omega_1$ has
nonempty interior, $\lambda_1^{T_j}(\Omega_1)>
\lambda_1^{T_j}(\Omega_2)$. Hence, we deduce the existence of
$$
\lambda_1^{T_j}(M) = \lim_{\mu \ra
+\infty}\lambda_1^{T_j}(\Omega_{\mu}),
$$
where $\{\Omega_\mu\}$ is any exhaustion of $M$ by means of
increasing, relatively compact domains with smooth boundary. The
next result is substantially an application of the result of
Moss-Piepenbrink \cite{mosspie}, slightly modified according to
Fischer-Colbrie and Schoen \cite{fishschoen} and Fischer-Colbrie
\cite{fish} (consult also \cite{vanishing}, Chapter $3$ and, for
the case of $L_1$, \cite{elbert}).
\begin{proposition}\label{firsteig}
Let $M$ be a Riemannian manifold and let $T_j$ be as above. The
following statements are equivalent:
\begin{itemize}
    \item[(i)] $\lambda_1^{T_j}(M) \geq 0$;
    \item[(ii)] there exists $u \in C^{\infty}(M)$, $u > 0$ solution of $T_j u=0$ on $M$.
\end{itemize}
Furthermore, there exists a compact set $K \subset M$ and $u \in
C^{\infty}(M \backslash K)$, $u >0$ solution of $T_ju =0$ on $M
\backslash K$ if and only if $\lambda_1^{T_j}(M \backslash K) \geq
0$.
\end{proposition}

Next, we shall need to consider the following Cauchy problem;
here, as usual, $\erre^{+} = (0, + \infty)$ and $\erre_0^{+} = [0,
+ \infty)$.
\begin{equation}\label{oscil}
\begin{cases}
(v(t)z'(t))'+ A(t)v(t)z(t) = 0 \qquad \text{on } \erre^{+} \\
z'(t) = \mathrm{O}(1) \ \ \mathrm{as}\ t \downarrow 0^+, \ \
z(0^+)=z_0 > 0
\end{cases}
\end{equation}
where $A(t)$ and $v(t)$ satisfy the following conditions:
\begin{itemize}
\item[(A1)] $A(t) \in L_{\text{loc}}^{\infty}(\erre_0^{+})$, $\
A(t) \geq 0$, $\ A\not\equiv 0 \ \text{in
L$_{\text{loc}}^{\infty}$ sense}$; \item[(V1)] $v(t) \in
L_{\text{loc}}^{\infty}(\erre_0^{+})$, $\ v(t) \geq 0$, $\
\dfrac{1}{v(t)} \in L_{\text{loc}}^{\infty}(\erre^{+})$;
\item[(V2)] there
exists $a \in \erre^+$ such that $v$ is increasing on $(0,a)$ and\\
$\lim_{t \ra 0^+} v(t) = 0$.
\end{itemize}

Observe that $(V2)$ has to be interpreted as there exists a
version of $v$ which is increasing near $0$ and whose limit as $ t
\ra 0^+$ is $0$.

By Proposition $A.1$ of \cite{bimari} under the above assumptions
\eqref{oscil} has a solution $z(t) \in
\text{Lip}_{\text{loc}}(\erre_0^+)$ (and condition $z'(t)=
\text{O}(1)$ as $t \downarrow 0^+$ is satisfied in an appropriate
sense). Furthermore by Proposition $A.3$ of \cite{bimari}, $z(t)$
has only isolated zeros. In case $1/v \in L^1((1,+\infty))$, by
Proposition $2.5$ of \cite{bimari} if, for some $T > 0$,
\begin{equation}\label{casoint}
\limsup_{t \ra \infty} \frac{\int_T^t \sqrt{A(s)} \di s}{-\frac 12
\log \int_t^{+\infty} \frac{\di s}{v(s)}} > 1
\end{equation}
then, every solution of
\begin{equation}
\begin{cases}
(v(t)z'(t))'+ A(t)v(t)z(t) = 0 \ \ \ \text{on} \ (t_0,+ \infty), \
t_0 > 0 \\
z(t_0)=z_0 > 0
\end{cases}
\end{equation}
has isolated zeros and is oscillatory. The same happens if
\begin{equation}\label{casononint}
\int^{+\infty} \frac{\di t}{v(t)} =+\infty \qquad \text{and}
\qquad \int^{+\infty}A(t)v(t)\di t =+\infty.
\end{equation}
(see Corollary $2.4$ of \cite{bimari}).

A final result that we shall use is the following computation.
(For a proof see \cite{ros} , \cite{alenccol}).

\begin{proposition}\label{jacobi}
Let $f:M\ra \erre^{m+1}$ be an isometric immersion of an oriented
hypersurface and $\nu:M\ra \esse^m$ its Gauss map. Fix $a \in
\esse^m$. Then
\begin{equation}
\begin{array}{l}
L_j \langle a , \nu \rangle = - \big(S_1 S_{j+1}-(j+2)S_{j+2}\big)
\langle a , \nu \rangle - \langle \nabla S_{j+1}, a \rangle;
\\[0.4cm]
L_j\langle f,\nu\rangle = -(j+1)S_{j+1} - \big(S_1
S_{j+1}-(j+2)S_{j+2}\big) \langle f , \nu \rangle - \langle \nabla
S_{j+1}, f \rangle.
\end{array}
\end{equation}
where $\langle, \rangle $ stands for the scalar product of vectors
in $\esse ^m \subset \erre^{m+1}$.
\end{proposition}
In particular, if $S_{j+1}$ is constant, we have $T_j\langle
a,\nu\rangle = 0$. Moreover, if $S_{j+1}\equiv 0$, $T_j\langle
f,\nu\rangle \equiv 0$.
\section{Proof of Theorem \ref{main}}
Fix an equator $E$ and reason by contradiction: assume that there
exists a sufficiently large geodesic ball $B_R$ such that, outside
$B_R$, $\nu$ does not meet $E$. In other words, $\nu (M \backslash
B_R)$ is contained in the open spherical caps determined by $E$.
Indicating with $a \in \esse^m$ one of the two focal points of
$E$, $\langle a, \nu(x) \rangle \neq 0$ on $M \backslash B_R $.

Let $\mathcal{C}$ be one of the (finitely many) connected
components of $M \backslash B_R$; then $\nu (\mathcal{C})$ is
contained in only one of the open spherical caps determined by
$E$. Up to replacing $a$ with $-a$, we can suppose $u= \langle a ,
\nu \rangle > 0$ on $\mathcal{C}$. Proceeding in the same way for
each connected component we can construct a positive function $u$
on $M \backslash B_R$. Since $S_{j+1}$ is constant, by Proposition
\ref{jacobi} we have that $u > 0$ satisfies
$$
T_j u = L_j u + \big(S_1 S_{j+1}-(j+2)S_{j+2}\big)u=0
$$
on $M \backslash B_R$. Thus, by Proposition \ref{firsteig},
$\lambda_1 ^{T_j}(M \backslash B_R) \geq 0$.
We shall now show that the assumptions of the theorem contradict
this fact. As already stressed, the existence of an elliptic point
forces both $H_j$ and $H_{j+1}$ to be positive. Fix a radius $0 <
R_0 < R$ and let $K_j$ be a smooth positive function on $M$ such
that
\begin{equation}
K_j(x) = \left\{ \begin{array}{ll}
1 & \textrm{on $B_{R_0/2}$}\\[0.2cm]
(m-j)S_j & \textrm{on $M \backslash B_{R_0}$}
\end{array} \right.
\end{equation}
Next, we define
\begin{equation}
v_j(t)=\int_{\partial B_t} K_j
\end{equation}
Using Proposition $1.2$ of \cite{bimari} we see that $v_j(t)$ satisfies
$(V1)$ with $v_j(t) > 0$ on $\erre^+$ and $(V2)$. Next, we define
\begin{equation}\label{defA}
A(t)=\frac{1}{v_j(t)}\int_{\partial B_t} S_1 S_{j+1}-(j+2)S_{j+2}.
\end{equation}
Then, repeated applications of Newton inequalities give
\begin{equation}\label{newtonineq}
H_1H_{j+1}-H_{j+2} \geq 0.
\end{equation}
Thus, using \eqref{newtonineq}
\begin{equation}
\begin{array}{l}\label{stimadasotto}
\displaystyle S_1S_{j+1}-(j+2)S_{j+2}= m {m \choose
{j+1}}H_1H_{j+1}-(j+2){ m
\choose {j+2} }H_{j+2}=\\[0.4cm]
\displaystyle ={ m \choose {j+1}} (mH_1 H_{j+1}-(m-j-1)H_{j+2}) \\[0.4cm]
\displaystyle \ge {m\choose
j+1}\left[m-\frac{m-j-1}{j+2}\right]H_1H_{j+1} = (j+1){m+1\choose
j+2} H_1H_{j+1}\ge 0.
\end{array}
\end{equation}
This implies $A(t)\ge 0$, and
$$
A(t)v_j(t) \ge (j+1){m+1\choose j+2}H_{j+1} \int_{\partial
B_t}H_1=(j+1){m+1\choose j+2}H_{j+1} v_1(t).
$$
If $1/v_j \not\in L^1((1,+\infty))$, then under \eqref{assunz}
$(i)$ and by the coarea formula we deduce $Av_j \not\in
L^1(\erre^+)$. Hence, we can apply \eqref{casononint} to deduce
that every solution of
\begin{equation}
\begin{cases}
(v_j(t)z'(t))'+ A(t)v_j(t)z(t) = 0 \ \ \ \text{on} \ (t_0,+
\infty), \
t_0 > 0 \\
z(t_0)=z_0 > 0
\end{cases}
\end{equation}
is oscillatory. The same conclusion holds when $1/v_j\in
L^1((1,+\infty))$. Indeed, from \eqref{defA}, \eqref{stimadasotto}
\begin{equation}\label{ine}
\frac{\int_T^t \sqrt{A(s)} \di s}{-\frac 12 \log
\int_t^{+\infty} \frac{\di s}{v_j(s)}} \ge 2\sqrt{(j+1){m+1\choose
j+2}H_{j+1}}\frac{\int_T^t \sqrt{\frac{v_1(s)}{v_j(s)}} \di
s}{-\log \int_t^{+\infty} \frac{\di s}{v_j(s)}}.
%
\end{equation}
Using De l'Hopital rule and \eqref{assunz} $(ii)$, \eqref{casoint}
is met. Let now $R<T_1<T_2$ be two consecutive zeros of $z(t)$
after $R$. Define
\begin{displaymath}
\psi(x) = \left\{ \begin{array}{ll} z(r(x)) & \text{on } \
\overline{B_{T_2}} \backslash
B_{T_1}\\[0.3cm]
0 & \text{outside } \ \overline{B_{T_2}} \backslash B_{T_1}.
\end{array} \right.
\end{displaymath}

Note that $\psi \equiv 0$ on $\partial (\overline{B_{T_2}}
\backslash B_{T_1})$, $\psi \in \text{Lip}_0 (M)$ and $\nabla
\psi(x)=z'(r(x)) \nabla r(x)$ where defined. Furthermore, by the
coarea formula and the definition of $A(t)$ we have
\begin{eqnarray*}
\int_M (S_1
S_{j+1}-(j+2)S_{j+2})\psi^2=\int_{T_1}^{T_2}z^2(t)\int_{\partial
B_t}(S_1 S_{j+1}-(j+2)S_{j+2}) \di t=\\
=\int_{T_1}^{T_2}z^2(t)A(t)v_j(t) \di t=(m-j)\int_M S_j A(r)
\psi^2.
\end{eqnarray*}
Thus, using \eqref{oscil}, the above identity and again the coarea
formula
\begin{eqnarray*}
&&\int_M \langle P_j(\nabla \psi),\nabla \psi
\rangle-(S_1S_{j+1}-(j+2)S_{j+2})\psi^2 \\
&& \leq \int_M
\mathrm{Tr}(P_j)|\nabla \psi|^2-(S_1S_{j+1}-(j+2)S_{j+2})\psi^2=\\
&&=\int_M (m-j)S_j|\nabla
\psi|^2-(S_1S_{j+1}-(j+2)S_{j+2})\psi^2\\
&&= (m-j)\int_{\overline{B_{T_2}} \backslash B_{T_1}} S_j
[(z')^2-A(t)z^2]=\\
&&=(m-j)\int_{T_1}^{T_2}[(z')^2-A(t)z^2]v_j(t) \di t =\\
&&=(m-j)\{z(t)z'(t)v_j(t)\big|_{T_1}^{T_2}-\int_{T_1}^{T_2}[(v_j(t)z'(t))'+A(t)v_j(t)z(t)]z(t)\di
t =0.
\end{eqnarray*}

It follows that
$$
\lambda_1^{T_j}(\overline{B_{T_2}}\backslash B_{T_1}) \leq
\frac{1}{\int_{M}\psi^2} \left\{\int_M \langle P_j(\nabla
\psi),\nabla \psi \rangle-(S_1S_{j+1}-(j+2)S_{j+2})\psi^2
\right\}=0.
$$
As a consequence $\lambda_1^{T_j}(M \backslash B_R)<0$, which
gives the desired contradiction.
\begin{remark}\label{remnonor}
\emph{As a matter of fact, the orientability of $M$ is not needed.
If $M$ is non orientable, $\nu$ is not globally defined. However,
changing the sign of $\nu$ does not change either the assumptions
or the conclusion of Theorem \ref{main}, since the antipodal map
on $\esse^m$ leaves each $E$ fixed. If $\langle a,\nu\rangle\neq
0$ on $M\backslash B_R$, the normal field $X=\langle a,\nu\rangle
\nu$ is nowhere vanishing and globally defined on $M\backslash
B_R$. This shows that, in any case, every connected component of
$M\backslash B_R$ is orientable.}
\end{remark}
\section{Proof of Theorem \ref{main2}}
Assume that, for some $K$, the tangent envelope of $M\backslash K$
does not coincide with $\erre^{m+1}$. By choosing cartesian
coordinates appropriately, we can assume
$$
0\not\in \bigcup_{p\in M\backslash K} T_pM.
$$
Then, the function $u=\langle f,\nu\rangle$ is nowhere vanishing
and smooth on $M\backslash K$. Up to changing the orientation,
$u>0$ on $M\backslash K$. By Proposition \ref{jacobi}, $T_ju =
-(j+1)S_{j+1}$=0. Note that here the assumption $H_{j+1}\equiv 0$
is essential. It follows that $\lambda_1^{T_j}(M\backslash K)\ge
0$. The rest of the proof is identical to that of Theorem
\ref{main}. Again, according to Remark \ref{remnonor} we can drop
the orientability assumption on $M$. Indeed, if the tangent
envelope of $M\backslash K$ does not cover $\erre^{m+1}$, the
vector field $X=\langle f,\nu\rangle \nu$ is a globally defined,
nowhere vanishing normal vector field on $M\backslash K$, hence
$M\backslash K$ is orientable.

\bibliographystyle{amsplain}
\bibliography{bibliojcurva}

\end{document}